\documentclass[english]{article}
\usepackage[T1, T2A]{fontenc}
\usepackage[russian]{babel}
\usepackage[utf8]{inputenc} 
\addto\captionsrussian{}
\usepackage{color}
\usepackage{amsmath,amssymb,amscd,mathrsfs}
\usepackage[T2A]{fontenc}
\usepackage[T1]{fontenc}
\usepackage{geometry}
\usepackage{graphicx} 
\usepackage{booktabs} 
\usepackage{array} 
\usepackage{paralist} 
\usepackage{verbatim} 
\usepackage{subfig} 
\usepackage{fancyhdr} 
\usepackage{bm}
\usepackage{sectsty}
\allsectionsfont{\sffamily\mdseries\upshape} 
\usepackage[nottoc,notlof,notlot]{tocbibind} 
\usepackage[titles,subfigure]{tocloft}

\usepackage[linewidth=1pt]{mdframed}
\usepackage[all]{xy}
\usepackage[normalem]{ulem}
\newtheorem{theorem}{Theorem}
\newtheorem{definition}{Definition}

\newtheorem{lemma}{Lemma}
\newtheorem{corollary}{Corollary}

\newenvironment{proof}{{\noindent\it Proof}\quad}{$\square$}

\title{Decomposing Dedekind Numbers: A Polynomial Representation with Powers of 2}
\author{Liu Yongqing }
\date{January 2023}

\begin{document}

\maketitle
\section{Introduction}
In 1897, Dedekind introduced the problem of counting the number of monotone Boolean functions, leading to the concept of Dedekind numbers \cite{dedekind1897zerlegungen}. This count, denoted by $\psi(n)$, grows rapidly with the number of variables $n$. Over the years, mathematicians like Kleitman, Gilbert, Hansel, Korshunov, Sapozhenko, and others explored formulas and estimates for Dedekind numbers \cite{gilbert1954lattice, Hansel1968, k1997, korshunov2003monotone, sapozhenko1989number, kleitman1969dedekind, kleitman1975dedekind}.

Kleitman established asymptotic bounds on $\psi(n)$ in 1969 \cite{kleitman1969dedekind}, and Korshunov provided explicit asymptotics in 1977 \cite{k1997}. Subsequent work by Sapozhenko in 1989 focused on antichains in partially ordered sets \cite{sapozhenko1989number}.

Efforts to compute exact Dedekind numbers beyond $n\leq 8$ faced challenges. Wiedemann computed the eighth Dedekind number in 1991 \cite{wiedemann1991computation}, but it took years for significant progress. Recent advances in 2023 by Pawelski \cite{pawelski2023divisibility}, Jäkel \cite{jakel2023computation}, and others \cite{van2023computation} yielded the ninth Dedekind number, showcasing the potential of computational optimizations.

Our contribution extends beyond previous studies. Our general Theorem 1 builds upon and generalizes existing formulas, including those presented in \cite{berman2021dedekind}. Introducing the concept of local Dedekind numbers, we provide a strategy for partitioning Dedekind numbers into finite terms, each comprising smaller Dedekind numbers. This partitioning allows us to decompose any Dedekind number into a polynomial of powers of 2.

In our ongoing research, we are actively exploring the determination of coefficients for these terms, with the goal of potentially deriving a more comprehensive formula for n-dimensional Dedekind numbers. 

\section{Definitions}

\begin{definition}[Cube]
    Let $E=\{0,1\}$, and define the order relation on it as $0\leq0, 0\leq 1, 1\leq 1$. 
    The set of all sequences of length $n$ formed by elements of $E$, denoted as $a=(a^1,\dots,a^n), a_i\in E, 0<i\leq n$, is denoted as $E^n$. The partial order on $E^n$ is defined as:
    $a\leq b$ if and only if $a^i\leq b^i$, $0<i\leq n$. We call $E^n$ an $n$-dimensional cube.
\end{definition}

\begin{definition}[Point Weight]
    The number of $1$s in an element $a\in E^n$ is called the weight of $a$, denoted as $|a|$.
\end{definition}

\begin{definition}[Comparable, Incomparable]
    Let $a,b\in E^n$, $a$ and $b$ are said to be comparable if $a\leq b$ or $b\leq a$. Otherwise, $a$ and $b$ are said to be incomparable.
\end{definition}

\begin{definition}[Cover]
    Let $a,b\in E^n$, $a$ is said to cover $b$ if $a\geq b$ and $|a|-|b|=1$, denoted as $b\prec a$.
\end{definition}

\begin{definition}[Antichain]
    Let $A$ be a poset, $B$ a subset of $A$. If any two elements in $B$ are comparable, then $B$ is called a chain. Conversely, if any two elements in $B$ are incomparable, then $B$ is called an antichain. A subset with only one element is both a chain and an antichain.
\end{definition}

\begin{definition}[Poset Dual]
    Given two posets $S = (X, \leq_1)$ and $S' = (X, \leq_2)$ with the same set of elements $X$, if for all $x, y \in X$, $x \leq_1 y$ if and only if $y \leq_2 x$, then the posets $S$ and $S'$ are said to be dual to each other.
\end{definition}

\begin{definition}[Monotone Boolean Function]
    A Boolean function with $n$ variables is a mapping $f:E^n\rightarrow E$. The function $f$ is monotone if for any $a,b\in E^n$, $a\leq b$ implies $f(a)\leq f(b)$. The set of all $n$-variable monotone Boolean functions is denoted by $M_n$.
\end{definition}

\begin{definition}[Dedekind Number]
    The number of monotone Boolean functions with $n$ variables is called the $n$th Dedekind number, denoted by $D_n$.
\end{definition}

\begin{definition}[Function Local Restriction]
    Let $S\subseteq E^n$ be a subset with the poset of $E^n$. Then the restriction of the $n$-variable Boolean function $f$ to $S$ is:
    \begin{equation*}
        f|_{S}: S\rightarrow E; x\mapsto f(x)
    \end{equation*}
    It is monotone if for any $a,b\in S$, $a\leq b$ implies $f|_{S}(a)\leq f|_{S}(b)$. For $n$-variable monotone Boolean functions $f, g$, if for all $x\in S$, $g|_{S}(x)=f|_{S}(x)$, then we say $g|_{S}=f|_{S}$.
    For all $n$-variable monotone Boolean functions, we call the number of all distinct monotone $f|_{S}$ the $n$-dimensional local Dedekind number on the subset $S$, denoted by $D(S)$, and its set is denoted by $M(S)$.
\end{definition} 

Note that $D_n = D(E^n)$.

\begin{definition}[Upper set, Lower set]
    Define the upper set of $a\in E^n$ as $S_{a,1}=\{b|b\geq a,b\in E^n\}$, 
    and the lower set as
    $S_{a,0}=\{b|b\leq a,b\in E^n\}$. 
\end{definition}

\begin{definition}[Subset generated by A, y]
    Let $A=\{a_1,\dots,a_k\}$ be a subset of $E^n$, $1\leq k\leq2^n$, 
    and $y=(y^1,\dots,y^k)$ be a $k$-tuple, where $y^j\in E$,
    define the subset of $E^n$ generated by $(A,y)$ as
    \begin{equation*}
        S_{A,y}=\bigcup_{a_i\in A} S_{a_i,y^i}
    \end{equation*}
\end{definition}

\section{Recursive Partition Formula of the Dedekind Number}

\begin{theorem}[Partition of the Dedekind number]
    For any subset $A=\{a_1,\dots,a_k\}$ of the poset $S\subseteq E^n$,
    the $n$-dimensional local Dedekind number $D(S)$ on $S$ can be decomposed into the sum of $D(A)$ parts:
    \begin{equation}
        D(S) = \sum_{f\in M(A)} D(S-S_{A,(f(a_1),\dots,f(a_k))})
    \end{equation}
    We call the set of all terms of the above polynomial the partition of $D(S)$ induced by $A$.
\end{theorem}

\begin{proof}
    This is equivalent to dividing the set $M(S)$ into several parts: given $f\in M(A)$, let $g$ be a monotonic Boolean function on $S$ such that $g(a_i)=f(a_i),1\leq i\leq k$. Due to monotonicity, the function values of $g$ on $S_{A,(f(a_1),\dots,f(a_k))}$ are determined, and the values on the remaining part $S-S_{A,(f(a_1),\dots,f(a_k))}$ are undetermined. Therefore, the number of $g$ that satisfy the condition is equal to the local Dedekind number of the remaining free part of $S$, that is, $D(S-S_{A,(f(a_1),\dots,f(a_k))})$.
\end{proof}

In $E^n$, let $|A|=1$ to obtain the following corollary, which has been previously proven by other methods in \cite{berman2021dedekind}.
\begin{corollary} 
    The $n$-dimensional Dedekind number $D_n=D(E^n-S_{a,1})+D(E^n-S_{a,0})$, where $a\in E^n$.
\end{corollary}

\begin{definition}[Complete Partition]
A subset $A\in S$ is said to completely partition $S$ if every element in the partition of $S$ induced by $A$ is either $D_k$ or a product of some $D_k$, where $k$ is an integer.
\end{definition}

Before delving into the concept of complete partition, we introduce the following relevant definitions, which will serve as the foundation for our subsequent derivations and proofs.

\begin{definition}[Cover-preserving Isomorphism]
Suppose $(P, \leq)$ and $(Q, \leq)$ are two posets. A mapping $f: P \rightarrow Q$ is called a cover-preserving isomorphism of posets if it satisfies the following two conditions:
\begin{enumerate}
    \item $f$ is a bijection.
    \item For all $x, y \in P$, if $x\prec y$ in $P$, then $f(x)\prec f(y)$ in $Q$. Conversely, if $f(x)\prec f(y)$ in $Q$, then $x\prec y$ in $P$.
\end{enumerate}
If a cover-preserving isomorphism exists between $P$ and $Q$, we denote this as $P \simeq Q$.
\end{definition}

\begin{definition}[The Triangle Subset $V_3$]
We define a special subset $A=\{a,b,v\}$ in $E^n$ such that $v\prec a, v\prec b$ and $|a|=|b|$, $a\ne b$. Its dual is $A'=\{a',b',v'\}$ such that $a'\prec v', b'\prec v'$ and $|a'|=|b'|$, $a'\ne b'$. We collectively refer to $A$ and $A'$ as $V_3$.
\end{definition}

\begin{theorem}[Conditions of Complete Partition]
A subset $A$ can completely partition $S\subseteq E^n$ if and only if there is no subset of $S-A$ that is cover-preserving isomorphic to $V_3$. 

If $S=E^n$, the condition is Equivalent to the following: subset $B$ of $S$ that is cover-preserving isomorphic to $E^2$, at least one of the following conditions is satisfied:
\begin{enumerate}
    \item At least one middle layer point of $B$ belongs to $A$.
    \item Both the top and bottom points of $B$ belong to $A$.
\end{enumerate}
\end{theorem}

\begin{proof}
Given a subset $A$ that completely partitions $S$, we proceed by contradiction. If there exists a subset $B$ that does not satisfy either condition, then $S$ contains a $V_3$. If a monotone function $f$ on the upper set of $V_3$ all takes the value $1$, on the lower set take the value $0$, and take $1$ on the points not comparable to $V_3$, then the number of this kind of monotone functions is equal to $D(V_3)$ which is not a $D_k$. Therefore, a subset $A$ that completely partitions $S$ must satisfy at least one of the two conditions. 

Conversely, if a subset $A$ satisfies at least one of the two conditions, then there is no $V_3$ in $S - A$. Therefore, for any $f\in M(S)$, and any subsets $B,B'$ of $S$ that is cover-preserving isomorphic to $E^2$, the function $f|_{B}$ will fix all points of $B$ or fix some points and leave subset $b$ that is cover-preserving isomorphic to $E\textsuperscript{1}, E\textsuperscript{0}$, and for $B'$, $f$ leaves $b'$.

When $b\simeq E^0$ and $b'\simeq E^0$, if they are comparable then $b\cup b'\simeq E^1$ else $b\cup b'\simeq E^0\cup E^0$, neither case violates the condition of complete partition. 
When $b\simeq E^0$ and $b'\simeq E^1$, they are not comparable, because if there is a pair pf elements $x\in b,x'\in b'$, such that $x\leq x'$ or $x\leq x'$, then there exists $y\in b'$ such that $V_3=\{x,x',y\}$, this contradicts the setting.
\end{proof}

\begin{lemma}[Minimal Complete Partition]
    Let $A$ be a subset of $E^n$, which completely partitions $E^n$. If there is no subgraph $B \subset A$ such that $B \simeq V_3$, then there is no $A' \subset E^n$ such that $A'$ conplete partition $E^n$ and $|A'| < |A|$.
\end{lemma}

\begin{proof}
Suppose $E^n$ is completely partitioned by $A$ and there is no subset in $(E^n - A)$ that is cover-preserving isomorphic to $V_3$. By the condition of complete partition, removing any point from $A$ would lead to the appearance of $V_3$ in $(E^n - A)$. Similarly, adding any point to $A$ would lead to the appearance of $V_3$ in $A$. So that adding a point $a\in (E^n-A)$ to $A$ cannot lead to the removal of any point $b \in A$, therefore the number of points in $A$ increases. Hence, a complete partition subset $A$ of $E^n$ with $V_3$ has more points than one without $V_3$. Therefore, the complete partition subset without $V_3$ is minimal.
\end{proof}

\begin{definition}[Upper and Lower Subcubes of $E^n$]
Let $E^n$ be an $n$-dimensional hypercube. For a given index $i, 0\leq i\leq n-1$, the set of points with the $i$-th coordinate equal to 1 forms the upper subcube, denoted as $E^n_{i,1}$. Similarly, the set of points with the $i$-th coordinate equal to 0 forms the lower subcube, denoted as $E^n_{i,0}$.
\end{definition}

\begin{theorem}[Construction of Minimal Complete Partition]
Let $E^n_{i,1}$ and $E^n_{i,0}$ for some $0\leq i\leq n-1$ be the upper subcube and the lower subcube of $E^n$ respectively. If $E^n_{i,1}$ (or $E^n_{i,0}$) can be completely partitioned by its subset $A_1$, and there is no $V_3$ in $A_1$, we take the complement $A_0$ of $A_1$ in the lower subcube $E^n_{i,0}$ (or $E^n_{i,1}$), then the subset $A_1\cup A_0$ can completely partition $E^n$.
\end{theorem}

\begin{proof}
Since there is no $V_3$ in $A_1$, there is no $V_3$ in the lower cube minus $A_0$, so the lower cube can be completely partitioned by $A_0$. For the edges connecting the upper and lower subcubes, each edge has one point belonging to $A_1 \cup A_0$, so no $V_3$ will appear. Therefore the union of $A_1$ and $A_0$ completely partitions $E^n$.
\end{proof}

\begin{lemma}[Size of the Complete Partition Subgraph]
    Suppose $A$ is a subset of $E^n$, and $A$ completely partitions $E^n$, then $|A|\geq 2^{(n-1)}$. If there exists a subgraph $B\simeq V_3$ in $A$, then $|A|>2^{(n-1)}$, otherwise $|A|=2^{(n-1)}$.
\end{lemma}

\begin{proof}
    The complete partition subgraph constructed in Theorem 3 has a size of $2^{(n-1)}$ and does not contain $V_3$, thus from lemma 1, we can conclude that if there exists a subgraph $B\simeq V_3$ in $A$, then $|A|>2^{(n-1)}$, otherwise $|A|=2^{(n-1)}$.
\end{proof}

In the subsequent discussion, we introduce a distinctive minimal partition subset that unveils a particular polynomial representation for Dedekind numbers. 

\begin{lemma}[A Specific Minimal Complete Partition]
Let $E^n$ be an $n$-dimensional hypercube with $n \geq 3$. When $n$ is even, its subsets $A=\{a\in E^n| |a|=0,2,\dots, n\}$ and $B=\{a\in E^n| |a|=1,3,\dots, n-1\}$ satisfy the construction of Theorem 3, and each forms a minimal complete partition subgraph of $E^n$. When $n$ is odd, its subsets $A'=\{a\in E^n| |a|=0,2,\dots, n-1\}$ and $B'=\{a\in E^n| |a|=1,3,\dots, n\}$ satisfy the construction of Theorem 3, and each forms a minimal complete partition subgraph of $E^n$.
\end{lemma}

\begin{proof}
In $(E^n - A)$, there are no pairs of points $a, b$ such that $a \prec b$. Therefore, in $(E^n - A)$, there are no subsets that are cover-preserving isomorphic to $V_3$. The same reasoning applies to $B$, $A'$, and $B'$.
\end{proof}

\begin{theorem}[Structure of the $D_n$ polynomial]
    Let $A$ be a subset of $E^n$ constructed by the method in Lemma3, which completely partitions $E^n$. Then in the polynomial
        \begin{equation}
            D(S) = \sum_{f\in M(A)} D(S-S_{A,(f(a_1),\dots,f(a_k))}),
        \end{equation}
    each term is a power of $E^0$.
\end{theorem}

\begin{proof}
    Consider a point $a$ in $(E^n - A)$. There exists $f \in M(A)$ such that $a \in (S-S_{A,(f(a_1),\dots,f(a_k))})$. Since $f$ is monotonic, for all such $b \in (E^n - A)$ with $b > a$, we have $f(b) = 1$, and for all such $c \in (E^n - A)$ with $c < a$, we have $f(c) = 0$. Therefore, both $b$ and $c$ do not belong to $(S-S_{A,(f(a_1),\dots,f(a_k))})$. As a result, $a$ is incomparable with all other points in $(S-S_{A,(f(a_1),\dots,f(a_k))})$. Therefore, it is impossible for a subset isomorphic to $E^k$ with $k \geq 1$ to exist in $(S-S_{A,(f(a_1),\dots,f(a_k))})$. However, $E^0$ that is incomparable with $a$ may appear. This completes the proof.
\end{proof}

\bibliographystyle{plain}
\bibliography{main}

\end{document}